%% file: JORDAN.TEX
\documentclass[12pt,a4paper,draft]{article}
\input{STYLE.TEX}
\title{On the theorem converse to Jordan's curve theorem.}
\author{Eugene Polulyakh}
\date{\makeatletter
        Institute of mathematics,\\
        National Acad. of Sci., Ukraine\\
        \medskip{\sl e-mail polulyah@imath.kiev.ua}
        \makeatother}
\begin{document}
\maketitle
\input{ABSTRACT.TEX}
\input{INTRO.TEX}
\input{SECT_1.TEX}
\input{SECT_2.TEX}
\input{SECT_3.TEX}

\input{REFREN.TEX}
\end{document}

%% file: STYLE.TEX
\usepackage{amssymb, amsfonts, amsmath}
\usepackage[mathscr]{eucal}
\usepackage{theorem}

\newcommand{\rr}{\mathbb R}
\newcommand{\nn}{\mathbb N}

\newcommand{\Cl}[1]{\overline{#1}}
\newcommand{\Front}{\mathscr F}
\newcommand{\Tau}{\mathcal T}
\newcommand{\LC}[1]{{\mathcal LC}(#1)}
\newcommand{\diam}{\mathop{\mathrm{diam}}}

\theoremstyle{plain}

\theorembodyfont{\sl}

\newtheorem{prop}{Proposition}[section]

\newtheorem{theorem}{Theorem}[section]

\newtheorem{lemma}{Lemma}[section]

\newtheorem{corr}{Corrolary}[section]

\theorembodyfont{\rm}

\newtheorem{examp}{Example}[section]

\newtheorem{defn}{Definition}[section]

\newtheorem{design}{Designation}[section]

\newtheorem{rem}{Remark}[section]

{\noindent%
{\bf Proof #1.} {} }%
{$\square$\medskip}

\newenvironment{proof_t}%
{\noindent%
$\bigtriangledown$ {}}%
{$\bigtriangleup$\medskip}

\newcounter{StepOfProof}
\renewcommand{\theStepOfProof}{\arabic{StepOfProof}}

\newenvironment{SteppedProof}[1][\hspace{-0.8ex}]%
{\noindent%
\setcounter{StepOfProof}{0}%
{\bf Proof #1}}%
{$\square$\medskip}

\newenvironment{StepOfProof}%
{\medskip
\noindent%
\refstepcounter{StepOfProof}%
{\bf \theStepOfProof}.}%
{}

%% file: ABSTRACT.TEX
\begin{abstract}
Theorem converse to Jordan's curve theorem says that {\it if a compact set $K$ has two
complementary domains in $\rr^{2}$, from each of which it is at every
point accessible, it is a simple closed curve}.

We show that the requirement of this theorem that {\it all} points of $K$
were accessible from {\it both} complementary domains is surplus and prove
one generalization of this theorem.
\end{abstract}

%% file: INTRO.TEX
\section*{Introduction.}

Approximately three years ago author was collided with the
following question: what are the sufficient conditions for a
compact set in the plane to be homeomorphic to the two-dimensional
disk? In addition, it was desirable to express supposed answer in
the terms of local properties of the frontier complementary
domains.

In this terms it is possible instead of considering a compact set
to pass to the investigation of its frontier. Then the previous
question could be restated in the form: what are the sufficient
conditions for a compact set in the plane to be a simple closed
curve?

Famous Jordan's curve theorem states that {\it a simple closed
curve in the plane has two complementary domains, of each of
which it is the complete frontier}. But the important conditions
given in this theorem are not sufficient. The counterexample is so
called polish curve (see the slight modification of it below,
example~\ref{example_2}).

Shoenflies theorem says that {\it an arbitrary homeomorphism of a
simple closed curve in the plane onto the unit circle could be
extended to a homeomorphism of the plane onto itself}. Every point
of the unit circle is accessible from the both complementary
domains. Hence an arbitrary point of a simple closed curve must be
accessible from the components of a complement. The conditions of
a sort are found to be the sufficient for a compact set in the
plane to be a Jordan curve.

So, the following theorem holds true: {\it if a closed set has two
complementary domains in $\rr^{2}$, from each of which it is at every
point accessible, it is a simple closed curve}. The author's
papers~\cite{Polulyakh} and~\cite{Polulyakh_2} was devoted to the proof of this theorem (it
was formulated there in rather different terms).

But the result given in~\cite{Polulyakh} and~\cite{Polulyakh_2} was turned out to be not
new. It is known as theorem converse to Jordan's curve theorem
(see~\cite{Kuratowski}, \cite{Newman}). Kuratovsky in his
book~\cite{Kuratowski} gives somewhat stronger statement: {\it if
for a closed set could be found two complementary domains in
$\rr^{2}$, from each of which it is at every point accessible, it
is a simple closed curve}. But the analysis of Newman's proof given
in~\cite{Newman} allows to conclude that it is suitable for the more
general case stated by Kuratovsky.

Author asks the natural question whether it is possible to weaken
the requirements of theorem converse to Jordan's curve theorem.

This problem appears to be nontrivial. For example the set of
points on the polish curve accessible from the both complementary
domains is open and dense (in particular it is massive). But this
set appears to be "badly disposed" in it.

Let $K$ be a compact set in the plane the complement of which has
two components $W_{1}$ and $W_{2}$. It proves out to be useful to
consider separately the sets $A_{1}$, $A_{2}$ of points of $K$,
accessible from $W_{1}$ and $W_{2}$, respectively. And if the both
of them are "well disposed" in $K$ then $K$ is a simple closed
curve.

It is interesting that if a set $K$ is connected and its
subsets $K \setminus A_{1}$ and $K \setminus A_{2}$ are
zero-dimensional then the sets $A_{1}$ and $A_{2}$ are "well
disposed" in $K$ and $K$ is a simple closed curve.

It appears that the concept of "well disposed" subset is connected
with one topological property of the space $K$ (with the topology
induced from the plane).

Namely, let $(X, \Tau)$ be a topological space. Denote by $\LC{X}$
the weakest from topologies on $X$ with the following property: for
each opened subset $W$ of $(X, \Tau)$ every connected component of
$W$ is opened in the topology $\LC{X}$.

In these terms the sets $A_{1}$ and $A_{2}$ will be "well disposed"
in $K$ if $A_{1}$ and $A_{2}$ are dense in $K$ in the topology
$\LC{K}$.

In what follows we will define the concept of "well disposed"
subset in a compact set in the plane and will give one
generalization of theorem converse to the Jordan's curve theorem.

Finally, I wish to express my deep gratitude to Yu.~B.~Zelinskiy
and V.~V.~Sharko for statement of a problem and formulation of
corrolary~\ref{corr_1}, and to M.~A.~Pankov for a number of
valuable remarks.

%% file: SECT_1.TEX
\section{Statement of main results.}

First we give some notations and definitions being used in what
follows.

Let $F$ be a set in the plane. We shall denote the frontier of $F$
by $\Front F$ and the closure of $F$ by $\Cl{F}$.

Let $x \in \rr^{2}$ and $\varepsilon > 0$. Denote
$$
U_{\varepsilon}(x) = \{ y \in \rr^{2} \,|\, \rho(x, y) <
\varepsilon \} \;.
$$

Let $U$ be an open subset of the plane. Three following definitions
(see~\cite{Newman}) are useful to study the local properties of
$\Front U$.

\begin{defn}
A simple continuous curve $\varphi : I \rightarrow \rr^{2}$ is
called an {\em end-cut of $U$ in a point $x \in \Front U$}, if
$\varphi(0) = x$ and $\varphi((0, 1]) \subset U$.
\end{defn}

\begin{defn}
Call a point $x \in \Front U$ {\em accessible from $U$}, if there
exists an end-cut of $U$ in the point $x$.
\end{defn}

\begin{defn}
{\em A cross cut of the domain $U$} is a simple continuous curve
$\psi : I \rightarrow \rr^{2}$ such that $\psi(0)$, $\psi(1) \in
\Front U$ and $\psi((0, 1)) \subset U$.
\end{defn}

And now we are ready to begin the main account.

\subsection{Definition and elementary properties of $d$-sets.}

Let $K$ be a compact subset of the plane dividing it into
two connected domains. Let us name such set
{\it two-sided}.

In what follows we shall denote complementary domains of a
two-sided set $K$ by $W_{1}$ and $W_{2}$. We shall also denote by
$A_{i}$ a set of points of $K$, accessible from $W_{i}$, $i = 1,
2$.

\begin{design}
Let $F$ be a set in $\rr^{2}$. Let $x \in F$ and $U$ be a
neighborhood of $x$ (not necessarily open). Designate
by $F(U, x)$ the connected component of $F \cap U$ which contains
$x$.
\end{design}

\begin{prop}\label{prop_1}
Let $F$ be a subset of the plane. Let $U_{1}$, $U_{2}$ be
two neighborhoods of point $x \in F$ and $U_{2} \subset U_{1}$.

Then $F(U_{2}, x) \subseteq F(U_{1}, x)$.
\end{prop}

\begin{proof_t}
The connected component of a point $x$ in $F \cap U_{1}$ is by
definition the maximal connected subset of $F \cap U_{1}$, which
contains $x$.

The connected set $F(U_{2}, x)$ lies in $F \cap U_{1}$ and
contains $x$. Therefore, $F(U_{2}, x) \subseteq F(U_{1}, x)$.
\end{proof_t}

\begin{defn}
Let $F$ be a set in the plane. We shall say that a subset $R$ of
$F$ is {\em sufficiently dense in $F$} if for all $x \in F$ and
every opened neighborhood $U$ of $x$ a set $F(U, x) \cap R$ is not
empty.
\end{defn}

\begin{rem}\label{rem_1}
It is easy to see that definition of sufficiently dense subset
could be reformulated in the following way: {\em a subset $R$ of a
set $F$ in the plane is sufficiently dense in $F$ if for any opened
set $U$ every connected component of a set $F \cap U$ contains a
point from $R$}.
\end{rem}

\begin{defn}\label{defn_1}
Two-sided set $K$ is called {\em $d$-set} if both $A_{1}$ and
$A_{2}$ are sufficiently dense in $K$.
\end{defn}

Let us deduce some simple properties of $d$-sets.

\begin{prop}\label{prop_2}
Let $K$ be a $d$-set.

For any $x \in K$ and its open neighborhood $U$ the sets
$A_{1}$ and $A_{2}$ are dense in $K(U, x)$.

Specifically, the sets $A_{1}$ and $A_{2}$ are dense in $K$.
\end{prop}

\begin{proof_t}
Let $y \in K (U, x)$ and $V$ be any open neighborhood
of $y$. Designate the set $V \cap U$ by $V_{0}$.

From property $d$ follows that there exist $y_{1}$, $y_{2}
\in K(V_{0}, y)$, such that $y_{i} \in A_{i}$, $i = 1, 2$.

But $K(V_{0}, y) \subseteq K(U, y) = K(U, x)$. Therefore,
$y_{1}$, $y_{2} \in K(U, x) \cap V$.
\end{proof_t}

\begin{prop}\label{prop_3}
Let $K$ be a two-sided set.

If the sets $A_{1}$ and $A_{2}$ are dense in $K$, then $K$ is
the common boundary of it's complementary domains, that is
\begin{equation}\label{eq_1}
K = \Front W_{1} = \Front W_{2} \;.
\end{equation}

In particular, any $d$-set complies with the relation~\ref{eq_1}.
\end{prop}

\begin{proof_t}
Obviously, $\Front W_{i} \subseteq K$, $i = 1, 2$.

On the other hand, $A_{i} \subseteq \Front W_{i} $, $i = 1, 2$. So
$K = \Cl{A_{i}} \subseteq \Front W_{i}$, $i = 1, 2$, since the
sets $\Front W_{i}$ are closed by the definition.
\end{proof_t}

For convenience of further account we shall give following

\begin{defn}
Let $K$ be a two-sided set. If $A_{1}$ and $A_{2}$ are dense in $K$
we shall say that $K$ is the {\em simple set}.
\end{defn}

It follows from proposition~\ref{prop_2} that $d$-sets are
simple.

\begin{prop}\label{prop_4}
Each simple set is connected.
\end{prop}

\begin{proof_t}
This statement is a direct corollary of the previous statement
and the following theorem: {\it if the domains $W_{1}$ and
$W_{2}$ in the plane do not meet, but $\Front W_{1} \subseteq
\Front W_{2}$, then $\Front W_{1}$ is connected}
(see~\cite{Newman}, theorem~V.14.1).
\end{proof_t}

\subsection{$d$-sets and theorem, inverse to Jordan's curve
theorem.}

The famous Jordan's curve theorem says, that {\it a simple closed
curve in the plane has two complementary domains, of each of
which it is the complete frontier}.

There is a natural question: under what conditions two-sided
set $K$ in the plane will be a Jordan curve?

Answer to this question gives the theorem, converse to
Jordan's curve theorem, which states that {\it a two-sided set,
all points of which are accessible from both components of
the complement, is a simple closed curve}
(see~\cite{Newman}, theorem~VI.16.1).

The requirement of this theorem that {\it all} points of $K$
were accessible from {\it both} complementary domains turns
out to be surplus.

How it is possible to weaken this condition?

Trivial direct check shows, that for any domain in $\rr^{2}$
points accessible from it are dense in the frontier of this
domain.  Therefore, taking account of proposition~\ref{prop_3}
it seems natural to require that sets $A_{i}$ of points,
accessible from the complementary domains $W_{i}$, $i =
1, 2$, were dense in $K$. This requirement could not be weakened,
as shows the following example.

\begin{examp}
$$
\begin{array}{lclc}
K & = & \{\mbox{the unit circle in } \rr^{2} \} & + \\
& + & \{\mbox{an isolated point in } W_{1} \} & + \\
& + & \{\mbox{an isolated point in } W_{2} \} \;. & \\
\end{array}
$$
\end{examp}

But even if the sets $A_{1}$ and $A_{2}$ are dense in
$K$ it is not sufficient for $K$ to be a Jordan curve.

\begin{examp}\label{example_2}
Consider a union $K$ of following sets
$$
\begin{array}{ccl}
K_{1} & = & \{ 0 \} \times [-1, 1] \;, \\
K_{2} & = & \left\{ (x, y) \in \rr^{2} \, | \, x \in (0, 1],
\, y = \sin \frac{\pi}{x} - \sin \pi x \right\}, \\
K_{3} & = & \left\{ (x, y) \in \rr^{2} \, | \, x \in (0, 1],
\, y = \sin \frac{\pi}{x} + \sin \pi x \right\} . \\
\end{array}
$$

It could be shown that $K$ is the simple set and each point of
the curve $K_{2} \cup K_{3}$ is accessible from both
components of complement.

However, the set $K$ is not arcwise connected (components
of it's linear connectivity are $K_{1}$ and $K_{2} \cup K_{3}$).

Also $K$ is not $d$-set. If $\varepsilon < 1$ then
$$
K(U_{\varepsilon}(0), 0) = \{ 0 \} \times
(-\varepsilon, \varepsilon)
$$
for circular neighborhood $U_{\varepsilon}(0)$,
and this set does not contain any point accessible from the
limited complementary domain of $K$.
\end{examp}

This example can be used as an argument for the benefit of
introduction of $d$-sets.

Let $F$ be a set in the plane. Since the topology on the plane has
a denumerable base of opened sets, we can find a denumerable subset
dense in $F$.

Hence, natural question arises: what is minimal cardinality of
sufficiently dense subsets of $F$?

The following example shows that there exists a simple set in
$\rr^{2}$ such that every its sufficiently dense subset is of
cardinality {\it continuum}.

\begin{examp} (See~\cite{Kuratowski}, section 48.I, examp. 4)
Let $C_{0}$ be the Cantor set situated at the axes $\xi_{1}$ of the
plane $(\xi_{1}, \xi_{2})$ and $C_{1}$ be the same set disposed on
the line $\xi_{2} = 1$. Connect every point of $C_{0}$ with the
appropriate point in $C_{1}$ by vertical interval. Add to $C_{0}$
the adjacent intervals of length $1/3$, $1/3^{3}, \ldots$ and to
$C_{1}$ add the adjacent intervals of length $1/3^{2}$, $1/3^{4},
\ldots$. We receive the continuum $K_{0}$ with a connected
complement.

\begin{center}
\input{PICTURE.TEX}
\end{center}

Let
$$
K_{1} = \left\{
(\xi_{1}, \xi_{2}) \in \rr^{2} \,|\,
(\xi_{1} - 1/2)^{2} + (\xi^{2} - 1)^{2} = 1/4 ,\; \xi_{2} \geq 1
\right\} \;.
$$
The arc $K_{1}$ is the cross-cut of $K_{0}$, so the continuum
$$
K = K_{0} \cup K_{1}
$$
is two-sided (see lemma~\ref{lemma_6}). It is not difficult to show
that the set $K$ complies with the equality~(\ref{eq_1}), hence $K$
is the simple set. In fact, $K$ is not a $d$-set.

It could be verified that for an opened set
$$
V = \rr \times (0, 1)
$$
every connected component of $V \cap K = C \times (0, 1)$ is $\{ c
\} \times (0, 1)$, where $c \in C$ is a certain element of the
Cantor set $C$.

So, the set $V \cap K$ has a {\it continuum} connected components
and from remark~\ref{rem_1} it follows that every sufficiently
dense subset of $K$ must have cardinality {\it continuum}.
\end{examp}

The main aim of the further considerations is to prove
following theorem and it's corollary.

\begin{theorem}\label{theorem_1}
Every $d$-set in the plane is a simple closed curve.
\end{theorem}

\begin{corr}\label{corr_1}
Let $K$ be a connected two-sided subset of the plane.
If sets $K \setminus A_{1}$ and $K \setminus A_{2}$ are
zero-dimensional then $K$ is a simple closed curve.
\end{corr}

\subsection{Proof of corollary~\ref{corr_1}.}

Corollary~\ref{corr_1} is a consequence of theorem~\ref{theorem_1}
and two following lemmas.

\begin{lemma}\label{lemma_A}
Let $F$ be a compact connected subset of the plane. Let $X$ be a
closed neighborhood of a point $x \in F$, such that $F \setminus
X \neq \emptyset$.

Then $F(X, x) \cap \Front X \neq \emptyset$.
\end{lemma}

\begin{proof_t}
Fix $y \in F \setminus X$.
It is known (see~\cite{Newman}, theorem~IV.5.1), that for
every $\varepsilon > 0$ points $x$ and $y$ can be
connected with the $\varepsilon$-chain in the set $F$,
i. e. final sequence
$$
x = z_{1}, z_{2}, \ldots, z_{k} = y
$$
of points could be found in $F$ to comply with the condition
$$
\rho (z_{i}, z_{i+1}) < \varepsilon \,, \quad i = 1, \ldots, k-1
\;.
$$

For every $n \in \nn $ select $1/n$-chain in $F$
$$
x_{1}^{n}, \ldots, x_{k(n)}^{n} \,, \quad n \in \nn \;,
$$
connecting points $x$ and $y$.

Since $y \notin X$, for every $n \in \nn$ in the
sequence $\{x_{i}^{n}\}_{i=1}^{k(n)}$ could be found an element
which does not belong to $X$. Therefore, for all $n \in \nn$
indexes $j(n) \in \{1, \ldots, k(n)-1\}$ are defined to
meet the conditions
\begin{itemize}
        \item [(i)] $x_{i}^{n} \in X$, $i = 1, \ldots, j(n)$;
        \item [(ii)] $x_{j(n)+1}^{n} \notin X$.
\end{itemize}
The set $F \cap X$ is compact, being the closed subset of the
compact set $F$. Therefore sequence $\{x_{j(n)}^{n}\}_{n \in
\nn}$ has a limit point $x_{0} \in X \cap F$. Without loss of
generality we shall suppose that $x_{j(n)}^{n} \rightarrow x_{0}$
(otherwise, it is always possible to pass to a subsequence).

Since $x_{j(n)+1}^{n} \notin X$ and $\rho(x_{j(n)}^{n},
x_{j(n)+1}^{n}) \rightarrow 0$ then $\rho(x_{j(n)}^{n}, \rr^{2}
\setminus X) \rightarrow 0$. Therefore, $\rho(x_{0}, \rr^{2}
\setminus X) = 0$ and $x_{0} \in \Front X$.

On the other hand, for every $\varepsilon > 0$ there exist $n_{1}
\in \nn$ and $n > n_{1}$, such that $1/n_{1} < \varepsilon$
and $\rho(x_{j(n)}^{n}, x_{0}) < \varepsilon$. Then points $x$
and $x_{0}$ could be connected in $F \cap X$ by the
$\varepsilon$-chain
$$
x = x_{1}^{n}, x_{2}^{n}, \ldots, x_{j(n)}^{n}, x_{0} \;.
$$
Consequently, (see~\cite{Newman}, theorem~IV.5.4) points
$x$ and $x_{0}$ belong to one connected component of the set $F
\cap X$.
\end{proof_t}

\begin{lemma}
Let $K$ be a connected two-sided subset of the plane.
If sets $B_{1} = K \setminus A_{1}$ and $B_{2} =
K \setminus A_{2}$ are zero-dimensional then $K$ is $d$-set.
\end{lemma}

\begin{proof_t}
Fix $x \in K$ and open neighborhood $U$ of $x$.
Let us show that $K(U, x) \cap A_{1} \neq \emptyset$.

The set $K$ divides the plane into two connected components.
Therefore it consists more than of one point. Take $\varepsilon >
0$ such that $U_{2 \varepsilon}(x) \subset U$ and $K \setminus
U_{2 \varepsilon}(x) \neq \emptyset$.

According to lemma~\ref{lemma_A} there exists $y \in
K(\Cl{U_{\varepsilon}(x)}, x) \cap \Front \Cl
{U_{\varepsilon}(x)}$. From proposition~\ref{prop_1} it follows
that $y \in K(U, x)$. Besides, $y \in K \setminus
U_{\varepsilon}(x)$.

The set $B_{1}$ is zero-dimensional. Therefore the space $B_{1}$
admits a base of open sets $\{V_{\alpha}\}_{\alpha \in \Lambda}$
such that $\Front V_{\alpha} = \emptyset$ in space $B_{1}$ (in
relative topology) for all $\alpha \in \Lambda$.

Let $x \in B_{1}$. Then $x \in V_{\beta}$ and $V_{\beta} \subset
U_{\varepsilon}(x)$ for some $\beta \in \Lambda$.

$V_{\beta}$ is an open-closed subset of space $B_{1}$.
If $V = B_{1} \setminus V_{\beta} \neq \emptyset$, the sets
$V$ and $V_{\beta}$ will form a partition of space $B_{1}$.

Assume that $K(U, x) \subset B_{1}$. Since $y \notin
U_{\varepsilon}(x)$ then $y \in V$ and $V \neq \emptyset$.
Therefore, the sets $V_{\beta} \cap K(U, x)$ and $V \cap
K(U, x)$ form a partition of the set $K(U, x)$ contrary to
it's connectedness.

So $K(U, x) \cap A_{1} \neq \emptyset$. The relation $K(U, x)
\cap A_{2} \neq \emptyset$ is proved similarly.

By virtue of arbitrariness in the choice of a point $x \in K$
and its neighborhood $U$ we conclude that $K$ is $d$-set.
\end{proof_t}

\subsection{On an accessibility of a point of a simple subset of
the plane from a component of its complement.}

It is easy to see that theorem~\ref{theorem_1} is
direct corollary of theorem, inverse to Jordan's curve
theorem, and the following statement.

\begin{theorem}\label{theorem_2}
Let $K$ be a simple subset of the plane. If for some
$x \in K$ and every open neighborhood $U$ of $x$
the set $K(U, x) \cap A_{1}$ is not empty, the point $x$ is
accessible from $W_{1}$.
\end{theorem}

\noindent
{\bf Plan of the proof of theorem~\ref{theorem_2}.}

If it is known that the point $x$ is accessible from $W_{1}$,
theorem holds true. Therefore we shall assume further that
$(K(U, x) \cap A_{1}) \setminus \{ x \} \neq \emptyset$ for
every open neighborhood $U$ of $x$.

Fix $x_{1} \in (K(U_{1}(x), x) \cap A_{1}) \setminus \{ x \}$.

Suppose for some $n \in \nn$ there are already selected $n$
different points $x_{1}, \ldots, x_{n}$ such that $x_{i} \in
K(U_{1/i}(x), x) \cap A_{1}$ and $x_{i} \neq x$, $i = 1, \ldots,
n$. Let
$$
\varepsilon_{1} = \min_{i \in \{1, \ldots, n\}} \rho(x_{i}, x) \;,
\quad \varepsilon = \min (\varepsilon_{1}, \frac{1}{n+1}) \;.
$$
Select $x_{n+1} \in (K(U_{\varepsilon}(x), x) \cap A_{1})
\setminus \{ x \}$. Proposition~\ref{prop_1} implies
the inclusion $x_{n+1} \in K(U_{1/(n+1)}(x), x)$.

By application of this procedure consequently for all $n \in
\nn$ we shall receive the sequence $\{x_{n}\}_{n \in \nn}$, with
pairwise different elements which complies with the condition
$$
x_{n} \in (K(U_{1/n}(x), x) \cap A_{1}) \setminus \{ x \} \,,
\quad n \in \nn \;.
$$
We fix for every $n \in \nn$ an end-cut $\alpha_{n} : I
\rightarrow \rr^{2}$ of the domain $W_{1}$ in the point $x_{n}$
(that is simple continuous curve, for which $\alpha_{n}(0) =
x_{n}$ and $\alpha_{n}((0, 1]) \subset W_{1}$).

Following statement which will be proved in the next section
holds true.

\begin{lemma}\label{lemma_1}
There exists a family of cross-cuts $\{\beta_{n} : I \rightarrow
\rr^{2}\}_{n=2}^{\infty}$ of the domain $W_{1}$ which complies
with the following properties
\begin{itemize}
        \item [(i)] a cross-cut $\beta_{n}$ connects points
                   $x_{n}$ and $x_{n+1}$, $n \geq 2$;
        \item [(ii)] $\beta_{n}(I) \subset U_{1/(n-1)}(x)$, $n
                   \geq 2$;
        \item [(iii)] for every $n \geq 2$ exists $\tau_{n} > 0$,
                   such that $\alpha_{2}([0, \tau_{2}])
                   \subset \beta_{2}(I)$ and $\alpha_{n}([0,
                   \tau_{n}]) \subset (\beta_{n-1}(I) \cap
                   \beta_{n}(I))$, $n \geq 3$.

\end{itemize}
\end{lemma}

Let $\{\beta_{n}\}_{n=2}^{\infty}$ and
$\{\tau_{n}\}_{n=2}^{\infty}$ be collections of cross-cuts and
parameter values from lemma~\ref{lemma_1}. From this lemma
follows that for every $n \geq 2$
$$
\alpha_{n}(\tau_{n}) \in U_{1/(n-1)}(x)
$$
and the points $\alpha_{n}(\tau_{n})$ and
$\alpha_{n+1}(\tau_{n+1})$ are contained in the same connected
component of $U_{1/(n-1)}(x) \cap W_{1}$. Therefore points
$\alpha_{n}(\tau_{n})$ and $\alpha_{n+1}(\tau_{n+1})$ can
be connected in the set $U_{1/(n-1)}(x) \cap W_{1}$ by a
simple polygonal line $J_{n}$ with the final number of links
(see~\cite{Newman}, theorem~V.6.3). We designate
$$
Q_{n} = \bigcup_{k=2}^{n} J_{k} \,, \quad
A = \bigcup_{n=2}^{\infty} Q_{n} = \bigcup_{n=2}^{\infty} J_{n}
\;.
$$
Obviously, the set $Q_{n}$ is connected and closed for every $n
\geq 2$. Therefore the set $A$ is connected as the union of the
connected sets with a common point $\alpha_{2}(\tau_{2})$.

Since $A \setminus Q_{n} \subset U_{1/(n-1)}(x)$, $n \geq
2$, then $\Cl{A} = A \cup \{ x \}$. Being a closure of the
connected set, $\Cl{A}$ is connected and on a construction
$\Cl{A} \setminus \{ x \} \subset W_{1}$.

The further proof is based on the following statement
(see~\cite{Newman}, theorem~VI.14.3).

\begin{lemma}\label{lemma_2}
Let $a \in \rr^{2}$ and $A$ be the union of a sequence of
segments, $x_{n} y_{n}$, such that $x_{n} \rightarrow a$, $y_{n}
\rightarrow a$. Then if $a$ is connected in $\Cl{A}$ to a point
$b$ ($b \neq a$), there is a simple arc in $\Cl{A}$ with
end-points $a$ and $b$.
\end{lemma}

From this lemma follows that there is a simple continuous curve
$\alpha : I \rightarrow \Cl{A} $ connecting points $x$ and
$\alpha_{2}(\tau_{2})$. Since $\Cl{A} \setminus \{ x \}
\subset W_{1}$, the arc $\alpha$ is an end-cut of the domain
$W_{1}$ in the point $x$.
$\square$
\medskip

So it suffices to prove lemma~\ref{lemma_1} for completion of
the proof of theorem~\ref{theorem_2}.

%% file: PICTURE.TEX
\begin{picture}(100, 50)

\thinlines

\put(0, 5){\line(1, 0){100}}
\put(5, 0){\line(0, 1){50}}

\thicklines

\put(5, 5){\line(0, 1){40}}
\put(35, 5){\line(0, 1){40}}
\put(65, 5){\line(0, 1){40}}
\put(95, 5){\line(0, 1){40}}
\put(35, 5){\line(1, 0){30}}

\multiput(0, 0)(60, 0){2}{%
\put(15, 5){\line(0, 1){40}}
\put(25, 5){\line(0, 1){40}}
\put(15, 45){\line(1, 0){10}}
}

\multiput(0, 0)(60, 0){2}{%
\multiput(0, 0)(20, 0){2}{%
\put(8.3, 5){\line(0, 1){40}}
\put(11.7, 5){\line(0, 1){40}}
\put(8.3, 5){\line(1, 0){3.4}}
}}

\multiput(0, 0)(60, 0){2}{%
\multiput(0, 0)(20, 0){2}{%
\multiput(0, 0)(6.6, 0){2}{%
\put(6.1, 5){\line(0, 1){40}}
\put(7.2, 5){\line(0, 1){40}}
\put(6.1, 45){\line(1, 0){1.1}}
}}}

\end{picture}

%% file: SECT_2.TEX
\section{Cross-cuts of complementary domains of simple sets in
the plane.}

\setcounter{equation}{0}

Lemma~\ref{lemma_1} from the previous section is based on one
local property of simple sets which will be studied now.

\begin{theorem}\label{theorem_3}
Let $K$ be a simple set, $x \in K$. Let for some
$\varepsilon > 0$ points $y_{1}$, $y_{2} \in K(\Cl{U_{\varepsilon}(x)}, x)$,
$y_{1} \neq y_{2}$, be accessible from a complementary domain
$W_{1}$ of $K$. Let
$$
\alpha_{i} : I \rightarrow \rr^{2}
$$
be an end-cut of $W_{1}$ in the point $y_{i}$, $i = 1, 2$.

For any $\varepsilon_{1} > \varepsilon$ there exists cross-cut
$$
\beta: I \rightarrow \rr^{2}
$$
of $W_{1}$ with end-points $y_{1}$, $y_{2}$, such that
\begin{itemize}
        \item [(i)] $\beta(I) \subset U_{\varepsilon_{1}}(x)$;
        \item [(ii)] $\alpha_{i}([0, \tau_{i}]) \subset \beta(I)$,
                $i = 1, 2$, for some $\tau_{1}$, $\tau_{2} \in
                (0, 1)$.
\end{itemize}
\end{theorem}

\begin{SteppedProof}[of theorem~\ref{theorem_3}] will be decomposed into
several steps.

Suppose $\varepsilon_{1} > \varepsilon$ is given. Fix
$\varepsilon_{2} \in (\varepsilon, \varepsilon_{1})$.

\begin{StepOfProof}
{\it There exists a cross-cut $\beta_{0} : I \rightarrow
\rr^{2}$ of $W_{1}$ which meets condition (ii) of theorem}.

Since $y_{1}$, $y_{2} \in U_{\varepsilon_{2}}(x)$, there
exists $\tau_{0} \in (0, 1)$ meeting the relations $\alpha_{i}([0,
\tau_{0}]) \subset U_{\varepsilon_{2}}(x)$, $i = 1, 2$, and
$\alpha_{1}([0, \tau_{0}]) \cap \alpha_{2} ([0, \tau_{0}]) =
\emptyset$.

Points $z_{1} = \alpha_{1}(\tau_{0})$, $z_{2} =
\alpha_{2}(\tau_{0})$ are contained in the open connected set $W_{1}$.
It is known that any domain in $\rr^{2}$ is arcwise connected.
Therefore, there exists simple arc
$$
\gamma : I \rightarrow W_{1}
$$
connecting points $z_{1}$ and $z_{2}$.

Designate
$$
\tau_{i} = \min \{ \tau \, | \, \alpha_{i}(\tau) \in \gamma (I) \}
\,, \quad i = 1, 2 \;.
$$
Nonempty closed sets $\{ 0 \} = \alpha_{i}^{-1}(y_{i})$ and
$\alpha_{i}^{-1}(\gamma (I))$ do not intersect,
hence $\tau_{i} \in (0, \tau_{0}]$, $i = 1, 2$.

Values $t_{i} \in [0, 1]$ are uniquely defined, such that
$\alpha_{i}(\tau_{i}) = \gamma (t_{i})$, $i = 1, 2$. Consider
a simple arc
$$
\beta_{0} : I \rightarrow \rr^{2} \;,
$$
$$
\beta_{0}(t) =
\left\{
        \begin{array}{ll}
             \alpha_{1}(4 t \tau_{1}) \,, & t \in [0, 1/4] \;, \\
             \gamma (2 (t-1/4) t_{2} + 2 (3/4-t) t_{1}) \,, &
                t \in [1/4, 3/4] \;, \\
             \alpha_{2}(4 (1-t) \tau_{2}) \,, & t \in [3/4, 1]
                \;. \\
        \end{array}
\right.
$$
This curve will be the desired cross-cut.

If, besides, $\beta_{0}(I) \subset U_{\varepsilon_{1}}(x)$,
then theorem is proved. Therefore, we shall suppose further that
$$
\beta_{0}(I) \setminus U_{\varepsilon_{1}}(x) \neq \emptyset \;.
$$
\end{StepOfProof}

\begin{StepOfProof}
We shall examine one special case. Let $K \subset
\Cl{U_{\varepsilon}(x)}$. It is clear that $K \subset
U_{\varepsilon_{1}}(x)$ then.

The following statement (see~\cite{Newman}, theorem~V.9.3)
holds true.

\begin{lemma}\label{lemma_4}
If a closed set $F \neq \rr^{2}$ is contained in a domain $D$ and
$D_{1}, D_{2}, \ldots$ are the components of $\rr^{2} \setminus F$,
the components of $D \setminus F$ are $D \cap D_{1}, D \cap D_{2},
\ldots$.
\end{lemma}

From this statement follows, that points $z_{1} =
\alpha_{1}(\tau_{0})$, $z_{2} = \alpha_{2}(\tau_{0})$ are contained
in the same connected component of the opened set
$U_{\varepsilon_{1}}(x) \setminus K$ (namely, in $W_{1} \cap
U_{\varepsilon_{1}}(x)$). Therefore, there exists a simple arc
$$
\gamma : I \rightarrow W_{1} \cap U_{\varepsilon_{1}}(x) \;,
$$
with the end-points $z_{1}$ and $z_{2}$.

Repeating argument of the previous step we conclude that there
exists a cross-cut
$$
\beta_{0} : I \rightarrow U_{\varepsilon_{1}}(x)
$$
of $W_{1} \cap U_{\varepsilon_{1}}(x)$ with the end-points
$y_{1}$, $y_{2}$, satisfying to a condition (ii) of theorem.

Hence, if $K \subset \Cl{U_{\varepsilon}(x)}$, then theorem is
true.

In all further argument we shall assume that
$$
K \setminus \Cl{U_{\varepsilon}(x)} \neq \emptyset \;.
$$
\end{StepOfProof}

\begin{StepOfProof}
{\it The set $K(\Cl{U_{\varepsilon}(x)}, x) \cup
\beta_{0}(I)$ is connected and divides the plane into two connected
domains $V_{1}$ and $V_{2}$, one of which is contained in $W_{1}$}.

First we prove following

\begin{lemma}\label{lemma_5}
Let $K$ be a simple set. For any proper closed subset
$\widetilde{K} \subsetneqq K$ it's complement $\rr^{2} \setminus
\widetilde{K}$ is connected.
\end{lemma}

\begin{proof_t}
Let $y \in K \setminus \widetilde{K}$. There exists $\varepsilon(y)
> 0 $, such that $U_{\varepsilon(y)}(y) \cap \widetilde{K} =
\emptyset$. According to proposition~\ref{prop_3}, we can find
$z_{1}$, $z_{2} \in U_{\varepsilon(y)}(y)$, such that $z_{i}(y) \in
W_{i}$, $i = 1, 2$.

Designate by $J_{i}(y)$ the segment with end-points $y$ and
$z_{i}(y)$, $i = 1, 2$. Obviously, $J_{i}(y) \cap \widetilde{K} =
\emptyset$, $i = 1, 2$. The set
$$
Q(y) = W_{1} \cup W_{2} \cup J_{1}(y) \cup J_{2}(y)
$$
is connected and $Q(y) \cap \widetilde{K} = \emptyset$.

We fix a point $z \in W_{1}$. Since $z \in Q(y)$ for every
$y \in K \setminus \widetilde{K}$, the set
$$
Q = \bigcup_{y \in K \setminus \widetilde{K}} Q(y) = \rr^{2}
\setminus \widetilde{K}
$$
is connected.
\end{proof_t}

Being the union of the connected sets which have a common point
$y_{1} = \beta_{0}(0)$, the set $K(\Cl{U_{\varepsilon}(x)}, x) \cup
\beta_{0}(I)$ is connected.

Under the definition $K(\Cl{U_{\varepsilon}(x)}, x)$ is the
connected component of a set $\Cl{U_{\varepsilon}(x)} \cap K$.
Therefore, it is closed in $\Cl{U_{\varepsilon}(x)} \cap K$.
Since $\Cl{U_{\varepsilon}(x)} \cap K$ is closed in $\rr^{2}$,
$K(\Cl{U_{\varepsilon}(x)}, x)$ is closed in $\rr^{2}$.

Applying lemma~\ref{lemma_5}, we conclude that the complement
$$
Q = \rr^{2} \setminus K(\Cl{U_{\varepsilon}(x)}, x)
$$
is connected.

Obviously, the curve $\beta_{0}$ is the cross-cut of $Q$.

The following statement holds true (see~\cite{Newman},
theorem~V.2.7).

\begin{lemma}\label{lemma_6}
If both the end-points of a cross-cut $L$ in a domain $D$ of
$\rr^{2}$ are on the same component of $\rr^{2} \setminus D$, $D
\setminus L$ has two components, and $L$ is contained in the
frontier of both.
\end{lemma}

This lemma implies that the set
$$
Q \setminus \beta_{0}(I) =
\rr^{2} \setminus (K(\Cl{U_{\varepsilon}(x)}, x) \cup
\beta_{0}(I))
$$
has two components $V_{1}$ and $V_{2}$, and that
\begin{equation}\label{eq_2}
\beta_{0}(I) \subset (\Front V_{1} \cap \Front V_{2}) \;.
\end{equation}

Since
\begin{equation}\label{eq_3}
W_{2} \cap \left(K(\Cl{U_{\varepsilon}(x)}, x) \cup \beta_{0}(I)
\right) = \emptyset \;,
\end{equation}
then $(V_{1} \cup V_{2}) \cap W_{2} \neq \emptyset$. Without loss
of generality, we can assume that $W_{2} \cap V_{2} \neq
\emptyset$.

$W_{2}$ is connected, therefore from the relation~(\ref{eq_3}) it
follows that $W_{2} \subset V_{2}$. Otherwise we should receive a
partition $W_{2} \cap V_{1}$, $W_{2} \cap V_{2}$ of $W _ {2} $.

Relation~(\ref{eq_1}) implies an equality $W_{1} = \rr^{2}
\setminus \Cl{W_{2}}$. On the other hand, $\Cl{W_{2}} \subset
\Cl{V_{2}} \subseteq \rr^{2} \setminus V_{1}$. Therefore,
$$
V_{1} \subseteq \rr^{2} \setminus \Cl{V_{2}} \subset
\rr^{2} \setminus \Cl{W_{2}} = W_{1} \;.
$$
\end{StepOfProof}

\begin{StepOfProof}
Now we shall describe a structure which will allow us to
construct a cross-cut $\beta$ of $Q \cap U_{\varepsilon_{1}}(x)$,
contained in a set $\beta_{0}(I) \cup V_{1}$.

It is known, that for any simple continuous curve in the plane
there exists a homeomorphism of $\rr^{2}$ onto itself, mapping this
curve onto a segment (see~\cite{Newman}, paragraph~VI.17).
Fix homeomorphism
$$
f_{0} : \rr^{2} \rightarrow \rr^{2}
$$
which maps the curve $\beta_{0}$ onto a segment $[-2, 2] \times \{
0 \}$.

Easy to see that it is possible to select such homeomorphism
$\varphi :  \rr \rightarrow \rr$, that the composition
$$
f = (\varphi \times id) \circ f_{0} : \rr^{2} \rightarrow \rr^{2}
$$
will satisfy to the following requirements:
\begin{itemize}
        \item[] $f(y_{1}) = f \circ \beta_{0}(0) = (-2 ,0)$,
                $f(y_{2}) = f \circ \beta_{0}(1) = (2 ,0)$;
        \item[] $f \circ \beta_{0}(1/4) =
                 f \circ \alpha_{1}(\tau_{1}) = (-1 ,0)$,
                $f \circ \beta_{0}(3/4) =
                 f \circ \alpha_{2}(\tau_{2}) = (1, 0)$.
\end{itemize}
Under these conditions it is obvious (see step 1) that
$$
\begin{array}{ccccccc}
f \circ \beta_{0}([0, 1/4]) & = & f \circ \alpha_{1}([0,
\tau_{1}]) & = & [-2, -1] \times \{ 0 \} & \subset &
f(U_{\varepsilon_{2}}(x)) \;, \\
f \circ \beta_{0}([3/4, 1]) & = & f \circ \alpha_{2}([0,
\tau_{2}]) & = & [1, 2] \times \{ 0 \} & \subset &
f(U_{\varepsilon_{2}}(x)) \;. \\
\end{array}
$$

Since $f$ is the homeomorphism and the set $U_{\varepsilon_{2}(x)}$
is open in $\rr^{2}$, there exists $\delta_{1} > 0$ meeting the
relations
$$
U_{\delta_{1}}(f \circ \beta_{0}(1/4)) \subset
f(U_{\varepsilon_{2}}(x)) \quad \mbox {and} \quad
U_{\delta_{1}}(f \circ \beta_{0}(3/4)) \subset
f(U_{\varepsilon_{2}}(x)) \;.
$$

Fix also $\delta_{2} > 0$, such that
$$
U_{\delta_{2}}([-1, 1] \times \{ 0 \}) \cap
f(K(\Cl{U_{\varepsilon}(x)}, x)) = \emptyset \;.
$$
We can do it, since $f(K(\Cl{U_{\varepsilon}(x)}, x))$ and $f \circ
\beta_{0}([1/4, 3/4])$ are disjoint compact sets.

Let $\delta = \min (\delta_{1}, \delta_{2})$.

Consider the following subsets of $U_{\delta}([-1, 1] \times
\{ 0 \})$.
$$
\begin{array}{ccc}
P_{1} & = & (\{ -1 \} \times [0, \delta/2]) \cup ([-1, 1]
\times \{ \delta/2 \}) \cup (\{ 1 \} \times [0, \delta/2]) \;,
\\
P_{2} & = & (\{ -1 \} \times [-\delta/2, 0]) \cup ([-1, 1]
\times \{ -\delta/2 \}) \cup (\{ 1 \} \times [-\delta/2, 0])
\;. \\
\end{array}
$$
Let us set on these two polygonal lines a parametrization $\mu_{i}
: I \rightarrow P_{i}$, $i = 1, 2$, converting them to simple
continuous curves, such that
$$
\begin{array}{c}
\mu_{1}(0) = \mu_{2}(0) = f \circ \beta_{0}(1/4) = (-1, 0) \;, \\
\mu_{1}(1) = \mu_{2}(1) = f \circ \beta_{0}(3/4) = (1, 0) \;. \\
\end{array}
$$

Consider in addition a segment $J \subset U_{\delta}((-1, 0))$
with the end-points $(-1 - \delta/2, 0) \in f \circ \alpha_{1}
([0, \tau_{1}])$ and $(-1, \delta/2) \in \mu_{1}(I)$.

Let $\mu : S^{1} \rightarrow \rr^{2}$ be a simple closed
curve, defined by a relation
$$
\mu(t) = \left\{
\begin{array}{ll}
\mu_{1}(2t) \,, & t \in [0, 1/2] \;, \\
\mu_{2}(2-2t) \,, & t \in [1/2, 1] \;. \\
\end{array}
\right.
$$
The curve $\mu$ divides the plane into two domains. Designate a
limited component of the complement $\rr^{2} \setminus \mu(S^{1})$
by $V_{in}$ and unlimited component by $V_{o}$. Obviously,
$$
V_{in} = (-1, 1) \times (-\delta/2, \delta/2) \subset
U_{\delta}([-1, 1] \times \{ 0 \}) \;.
$$
Therefore, $f(K(\Cl{U_{\varepsilon}(x)}, x)) \subset V_{0}$.

Segment $[-1, 1] \times \{ 0 \} = f \circ \beta_{0}([1/4, 3/4])$
divides a rectangle $V_{in}$ into two domains. Designate by
$V_{in}^{+}$ the domain contained in the upper half-plane and the
other domain by $V_{in}^{-}$.

By a construction $V_{in} \cap f(K(\Cl{U_{\varepsilon}(x)},
x)) = \emptyset$. Therefore, according to relation~(\ref{eq_2}), $V_{in}
\cap f(\Front V_{1}) = V_{in} \cap f(\Front V_{2}) = f \circ
\beta((1/4, 3/4)) = (-1, 1) \times \{ 0 \} $. Consequently, one
from domains $V_{in}^{+}$, $V_{in}^{-}$ is contained in
$f(V_{1})$, the other is a subset of $f(V_{2})$.

We shall suppose that $V_{in}^{+} \in f(V_{1})$. Otherwise
we can replace $f$ by $i \circ f$, where $i : \rr^{2} \rightarrow
\rr^{2}$, $i : (x_{1}, x_{2}) \mapsto (x_{1}, -x_{2})$,
is an inversion relative to the first coordinate axes.

Since
$$
\mu_{i}(I) \cap f \left( K(\Cl{U_{\varepsilon}(x)}, x) \cup
\beta_{0}(I) \right) = \mu_{i}(0) \cup \mu_{i}(1) \,, \quad
i = 1, 2 \;,
$$
then
$$
\mu_{1}((0, 1)) \subset f(V_{1}) \,, \quad
\mu_{2}((0, 1)) \subset f(V_{2}) \;.
$$
\end{StepOfProof}

\begin{StepOfProof}
Consider a simple closed curve
$$
\gamma : S^{1} \rightarrow \rr^{2} \;,
$$
$$
\gamma (t) = x + \varepsilon_{2} \exp (2 \pi i t) \;,
$$
being the boundary of a neighborhood $U_{\varepsilon_{2}}(x)$.

Let us designate
$$
\theta = f^{-1} \circ \mu : \rr^{2} \rightarrow \rr^{2} \;.
$$

The following statement (see~\cite{Newman}, example~VI.16)
holds true.

\begin{lemma} [Ker\'{e}kj\'{a}rt\'{o}'s Theorem]\label{lemma_7}
If two simple closed curves, $J_{1}$ and $J_{2}$, have more than
one common point, all the residual domains of $J_{1} \cup J_{2}$
are Jordan domains.
\end{lemma}

Show that the curves $\gamma$ and $\theta$ meet the condition of
this lemma.

Since $\beta_{0}(I) \setminus U_{\varepsilon_{1}}(x) \neq
\emptyset$ (see step 1) and $\varepsilon_{2} < \varepsilon_{1}$
then $\beta_{0}(I) \setminus \Cl{U_{\varepsilon_{2}}(x)} \neq
\emptyset$. On a construction, $\beta_{0}(I) \setminus
\Cl{U_{\varepsilon_{2}}(x)} \subset \beta_{0}((1/4, 3/4)) \subset
f^{-1}(V_{in})$.

Fix $z_{0} \in \beta_{0}(I) \setminus \Cl{U_{\varepsilon_{2}}(x)}$.
The boundary of the domain $f^{-1}(V_{in})$ is a Jordan curve
$\theta$. Therefore $\theta(S^{1}) \setminus
\Cl{U_{\varepsilon_{2}}(x)} \neq \emptyset$. Really, the set
$f^{-1}(V_{o})$ is not limited, hence there exists $z_{0}' \in
f^{-1}(V_{o}) \cap (\rr^{2} \setminus
\Cl{U_{\varepsilon_{2}}(x)})$. A set $\rr^{2} \setminus
\Cl{U_{\varepsilon_{2}}(x)}$ is arcwise connected, so the simple
continuous curve $\psi : I \rightarrow \rr^{2} \setminus
\Cl{U_{\varepsilon_{2}}(x)}$ could be found to connect points
$z_{0}$ and $z_{0}'$. Domains $f^{-1}(V_{o})$ and $f^{-1}(V_{in})$
are disjoint, $z_{0} \in f^{-1}(V_{in})$, $z_{0}' \in
f^{-1}(V_{o})$, and $\theta(S^{1}) = \Front f^{-1}(V_{in}) = \Front
f^{-1}(V_{o})$. Therefore, $\psi(I) \cap \theta(S^{1}) \neq
\emptyset$.

Let $z_{1}$ be contained in $\psi(I) \cap \theta(S^{1}) \subset
\theta(S^{1}) \setminus \Cl{U_{\varepsilon_{2}}(x)}$.

On the other hand, $z_{2} = \beta_{0}(1/4) = \theta(0) \in
U_{\varepsilon_{2}}(x)$ on a construction.

Points $z_{1}$ and $z_{2}$ divide Jordan curve $\theta$ into
two simple arcs with common end-points
$$
\theta_{i} : I \rightarrow \rr^{2} \,, \quad i = 1, 2 \;,
$$
such that
$$
\begin{array}{c}
\theta_{i}(0) = z_{1} \in U_{\varepsilon_{2}}(x) \,, \quad
\theta_{i}(1) = z_{2} \in \rr^{2} \setminus
\Cl{U_{\varepsilon_{2}}(x)} \,, \quad i = 1, 2 \;; \\
\theta_{1}(I) \cap \theta_{2}(I) =
\{ z_{1} \} \cup \{ z_{2} \} \;, \\
\end{array}
$$
It is not difficult to see, that $\theta_{i} \cap \Front
U_{\varepsilon_{2}}(x) \neq \emptyset$, $i = 1, 2$. Consequently,
the intersection $\gamma(S^{1}) \cap \theta(S^{1})$ contains not
less than two points and for the curves $\gamma$ and $\theta =
f^{-1} \circ \mu$ lemma~\ref{lemma_7} is applicable.
\end{StepOfProof}

\begin{StepOfProof}
Let $V_{K}$ be a complementary domain of $\theta(S^{1})
\cup \gamma(S^{1})$ containing $K(\Cl{U_{\varepsilon}(x)}, x)$.

Let
$$
\eta : S^{1} \rightarrow \rr^{2}
$$
be a simple closed curve bounding the domain $V_{K}$.

On the construction a connected set
$$
\begin{array}{cccc}
K_{0} & = &
K(\Cl{U_{\varepsilon}(x)}, x) \cup
\beta_{0}([0, 1/4) \cup (3/4, 1]) & = \\
& = & K(\Cl{U_{\varepsilon}(x)}, x) \cup
\alpha_{1}([0, \tau_{1})) \cup \alpha_{2}([0, \tau_{2})) & \\
\end{array}
$$
is contained in $V_{K}$. Therefore, $\beta_{0}(1/4)$,
$\beta_{0}(3/4) \in \Front V_{K} = \eta(S^{1})$. Similarly,
$$
(K_{0} \cup f^{-1}(J)) \setminus f^{-1}((-1, \delta/2))
$$
is the subset of $V_{K}$, hence $x_{0} = f^{-1}((-1, \delta/2)) \in
\Front V_{K} = \eta(S^{1})$. Since $x_{0} \in f^{-1} \circ
\mu_{1}((0, 1)) \subset V_{1}$, then $V_{1} \cap \Front V_{K} \neq
\emptyset$.

The curve $\eta$ is divided by points $\beta_{0}(1/4)$ and $\beta_{0}(3/4)$
into two arcs
$$
\eta_{i} : I \rightarrow \rr^{2} \,, \quad i = 1, 2 \;,
$$
such that
$$
\begin{array}{c}
\eta_{i}(0) = \beta_{0}(1/4) \,, \quad \eta_{i}(1) =
\beta_{0}(3/4) \,, \quad i = 1, 2 \;, \\
\eta_{1}(I) \cap \eta_{2}(I) =
\{ \beta_{0} (1/4) \} \cup \{ \beta_{0}(3/4) \} \;. \\
\end{array}
$$

Let $\eta_{1}$ be that from these two curves, which contains the
point $x_{0}$.
\end{StepOfProof}

\begin{StepOfProof}
We shall show, that $\eta_{1}((0, 1)) \subset V_{1}$.

The set $V_{1}$ is connected, so it will be enough to prove that
$\eta_{1}((0, 1)) \cap \Front V_{1} = \emptyset$. We shall check a
somewhat stronger relation
$$
\eta_{1}((0, 1)) \cap \left( K(\Cl{U_{\varepsilon}(x)})
\cup \beta_{0}(I) \right) = \emptyset \;.
$$

Since $K_{0} \subset V_{K}$, then $\eta_{1}(I) \cap K_{0} =
\emptyset$ and
$$
\eta_{1}((0, 1)) \cap \left( K(\Cl{U_{\varepsilon}(x)})
\cup \beta_{0}([0, 1/4] \cup [3/4, 1]) \right) = \emptyset \;.
$$

On the other hand, the sets $K(\Cl{U_{\varepsilon}(x)}, x)$ and
$\beta_{0}((1/4, 3/4))$ are contained in the different components
of $\rr^{2} \setminus \theta(S^{1})$, namely,
$$
\begin{array}{c}
K(\Cl{U_{\varepsilon}(x)}, x) \subset V_{K}
\subset f^{-1} (V_{o}) \;, \\
\beta_{0}((1/4, 3/4)) \subset f^{-1}(V_{in}) \;. \\
\end{array}
$$
Therefore,
$$
\eta_{1}((0, 1)) \cap \beta_{0}((1/4, 3/4)) \subset
\Cl{V_{K}} \cap f^{-1}(V_{in}) \subset
\Cl{f^{-1}(V_{o})} \cap f^{-1}(V_{in}) = \emptyset \;.
$$
\end{StepOfProof}

\begin{StepOfProof}
Recall that $V_{1} \subset W_{1}$ and arcs $\alpha_{1}$,
$\alpha_{2}$ are end-cuts of the domain $W_{1}$ in the points
$y_{1}$ and $y_{2}$, respectively. Since
$$
\begin{array}{c}
\alpha_{i}([0, \tau_{i})) \subset V_{K} \subset
U_{\varepsilon_{2}}(x) \,, \quad i = 1, 2 \;, \\
\eta_{1}(I) \subset \Front V_{K} \;, \\
\alpha_{1}(\tau_{1}) = \eta_{1}(0) \,, \quad
\alpha_{2}(\tau_{2}) = \eta_{1}(1) \;, \\
\end{array}
$$
then a simple arc
$$
\beta : I \rightarrow \rr^{2} \;,
$$
$$
\beta(t) =
\left\{
       \begin{array}{ll}
             \alpha_{1}(4 t \tau_{1}) = \beta_{0}(t) \,, &
                t \in [0, 1/4] \;, \\
             \eta(t/2+1/4) \,, & t \in [1/4, 3/4] \;, \\
             \alpha_{2}(4 (1-t) \tau_{2}) = \beta_{0}(t) \,, &
                t \in [3/4, 1] \\
       \end{array}
\right.
$$
is the cross-cut of the domain $W_{1}$, which complies with all
the requirements of theorem~\ref{theorem_3}.
\end{StepOfProof}
\end{SteppedProof}

\noindent
{\bf Proof of lemma~\ref{lemma_1}.}
From proposition~\ref{prop_1} follows, that
$$
K(U_{1/(n+1)}(x), x) \subset K(U_{1/n}(x), x) \subset
K(\Cl{U_{1/n}(x)}, x)
$$
for all $n \in \nn $. Therefore $x_{n}$, $x_{n+1} \in
K(\Cl{U_{1/n}(x)}, x)$, $n \in \nn$.

Apply now for all $n \geq 2$ theorem~\ref{theorem_3} to the values
$\varepsilon = 1/n$, $\varepsilon_{1} = 1/(n-1)$, points $y_{1} =
x_{n}$, $y_{2} = x_{n+1}$ and end-cuts $\alpha_{n}$,
$\alpha_{n+1}$. We receive a sequence $\{ \beta_{n} \} _
{n=2}^{\infty}$ of cross-cuts of the domain $W_{1}$ and two
sequences of positive numbers $\{ \tau_{n}' \}_{n=2}^{\infty}$, $\{
\tau_{n}''\}_{n=2}^{\infty}$, connected by relations
\begin{itemize}
        \item [(i)] $\beta_{n}(I) \subset U_{1/(n-1)}(x)$, $n \geq 2$;
        \item [(ii)] $(\alpha_{n}([0, \tau_{n}']) \cup
                \alpha_{n+1}([0, \tau_{n}''])) \subset
                \beta_{n}(I)$, $n \geq 2$.
\end{itemize}
Assume $\tau_{2} = \tau_{2}'$,
$\tau_{n} = \min(\tau_{n}', \tau_{n-1}'')$, $n \geq 3$.
Then the sequences $\{ \beta_{n} \}$ and $\{\tau_{n}\}$
satisfy to lemma~\ref{lemma_1}.
$\bigtriangleup$

%% file: SECT_3.TEX
\section{Concluding remarks.}

Theorems~\ref{theorem_1} and~\ref{theorem_2} could be
reformulated in an other way. In order to do that we shall
consider different topology on a two-sided subset of the plane.

\begin{defn}\label{defn_3_1}
Let $(X, \Tau)$ be a topological space. Let $\{ U_{\alpha}
\}_{\alpha \in A}$ be a base of opened sets in $X$.

For every $\alpha \in A$ let $\{V_{\alpha \beta}\}_{\beta \in
B(\alpha)}$ be the family of all connected components of the set
$U_{\alpha}$.

The topology induced by the family $\{ V_{\alpha \beta}
\}_{\alpha \in A,\, \beta \in B(\alpha)}$ will be designated by
$\LC{X, \Tau}$ (or by $\LC{X}$ if the original topology on
$X$ is clear from context).
\end{defn}

\begin{rem}
This topology could be informally characterized by the following
property: for each opened subset $W$ of $(X, \Tau)$ every
connected component of $W$ is opened in the topology $\LC{X,
\Tau}$.

Moreover, as it could be seen from definition~\ref{defn_3_1}
$\LC{X, \Tau}$ is the weakest topology on $X$ meeting this
property.

So, this remark could be used as the alternative definition of
the topology $\LC{X, \Tau}$.
\end{rem}

\begin{rem}
It could be easily verified that the space $(X, \Tau)$ is locally
connected if and only if the topologies $\Tau$ and $\LC{X,
\Tau}$ coincide.
\end{rem}

Using definition~\ref{defn_3_1} theorems~\ref{theorem_1}
and~\ref{theorem_2} could be reformulated in the following form

\begin{theorem}
Let $K$ be a two-sided subset of the plane. If the sets $A_{1}$
and $A_{2}$ are dense in $K$ in the topology $\LC{K}$ then $K$
is the simple closed curve.
\end{theorem}

\begin{theorem}
Let $K$ be a simple subset of the plane. Then the sets $A_{1}$
and $A_{2}$ are closed in the topology $\LC{K}$.
\end{theorem}

For a special class of spaces the topology
$\LC{X, \Tau}$ is mentioned in~\cite{Kuratowski}, chapter~6,
section~49.VII (see also~\cite{Mazurkevich}).

Namely, let $(X, \rho)$ be a metric space with the following
property: for every $a, b \in X$ there exists a connected subset $A
\subset X$ such that $a, b \in A$ and
$$
\diam A = \sup_{x, y \in A} \rho (x, y) < \infty \;.
$$

Then so called {\sl relative distance}
$$
\rho_{r} : X \times X \rightarrow \rr_{+} \;,
$$
$$
\rho_{r} (a, b) = \inf \{ \diam A \;|\; a, b \in A \mbox{ and }
A \mbox{ is connected}\}
$$
is well defined and appears to be a metric.

\begin{prop}
Topology induced by the distance function $\rho_{r}$ coincides
with the topology $\LC{X}$.
\end{prop}

\begin{proof_t}
Fix $x \in X$. For $\varepsilon > 0$ designate by
$X(U_{\varepsilon}(x), x)$ the connected component of $x$ in the
set $U_{\varepsilon}(x)$.

Then in accord with definition~\ref{defn_3_1} the family $\{
X(U_{\varepsilon}(x), x) \}_{\varepsilon > 0}$ forms a base of
opened neighborhoods of the point $x$ in the topology $\LC{X}$.

On the other hand let us denote
$$
V^{r}_{\varepsilon}(x) = \{ y \in X \;|\; \rho_{r}(x, y) <
\varepsilon \}
$$
for every $\varepsilon > 0$.

Then by easy direct verification we receive
$$
V^{r}_{\varepsilon/3}(x) \subseteq X(U_{\varepsilon}(x), x)
\subseteq V^{r}_{3\varepsilon}(x)
$$
for every $\varepsilon > 0$.

Q. E. D.
\end{proof_t}